\documentclass[12pt,twoside]{article}
\usepackage{latexsym, amssymb,amsthm,amsmath,xypic,lscape,epsfig,setspace,tikz,slashed}
\usepackage[retainorgcmds]{IEEEtrantools}
\usetikzlibrary{decorations.markings}
 \usetikzlibrary{arrows}
\usepackage[pdftex]{hyperref}
\usepackage[margin=10pt,font=small,labelfont=bf]{caption}

\theoremstyle{plain}
\newtheorem{theorem}{Theorem}[section]
\newtheorem{proposition}[theorem]{Proposition}
\newtheorem{lemma}[theorem]{Lemma}
\newtheorem{question}[theorem]{Question}

\theoremstyle{definition}
\newtheorem{definition}[theorem]{Definition}

\newtheorem{example}[theorem]{Example}

\newenvironment{renumerate}%
{%
\begin{enumerate}}%
{\end{enumerate}%
}%

\newenvironment{remark}%
{\vskip6pt%
\noindent%
{\it Remark.}}%
{\vskip6pt}

{\vskip6pt%
\noindent%
{\it Remarks}. %
\begin{renumerate}}%
{\end{renumerate}\vskip6pt}

\makeatletter
\def\Ddots{\mathinner{\mkern1mu\raise\p@
\vbox{\kern7\p@\hbox{.}}\mkern2mu
\raise4\p@\hbox{.}\mkern2mu\raise7\p@\hbox{.}\mkern1mu}}
\makeatother

\newcommand{\R}{\text{${\mathbb R}$}}
\newcommand{\C}{\text{$\mathbb C$}}
\newcommand{\NN}{\text{$\mathbb N$}}

\newcommand{\Z}{\text{$\mathbb Z$}}

\renewcommand{\tilde}{\widetilde}

\renewcommand{\L}{\text{$\mathcal{L}$}}

\newcommand{\gf}{\text{$\varphi$}}

\newcommand{\del}{\text{$\partial$}}

\newcommand{\tensor}{\otimes}

\newcommand{\mc}[1]{\text{$\mathcal{#1}$}}

\newcommand{\noqed}{\let\qed\relax}

\newcommand{\IP}[1]{\langle #1 \rangle}

\newcommand{\Cour}[1]{[\![#1]\!]}

\date{} \usepackage{color} \definecolor{tocolor}{rgb}{.1,.1,.5}
\definecolor{urlcolor}{rgb}{.2,.2,.6}
\definecolor{linkcolor}{rgb}{.1,.1,.6}
\definecolor{citecolor}{rgb}{.6,.2,.1}
\hypersetup{backref=true, colorlinks=true, urlcolor=urlcolor,
  linkcolor=linkcolor, citecolor=citecolor}

\begingroup
\hypersetup{linkcolor=blue}

\endgroup

\numberwithin{equation}{section}

\begin{document}

\title{Spherical T-duality\\
{\large Dimension change from higher degree $H$-flux}}
\author{Gil R. Cavalcanti\thanks{{\tt gil.cavalcanti@gmail.com}, Department of Mathematics,
Utrecht University, The Netherlands}, Bart Heemskerk\thanks{{\tt b.a.j.heemskerk@students.uu.nl}, Department of Mathematics,
Utrecht University, The Netherlands} and Bernardo Uribe\thanks{{\tt bjongbloed@uninorte.edu.co}, Mathematics and Statistics Department,
Universidad del Norte, Colombia}}
\maketitle

\abstract{Topological Spherical T-duality was introduced by Bouwknegt, Evslin and Mathai in \cite{BEM15} as an extension of topological T-duality from $S^1$-bundles to $\mathrm{SU}(2)$-bundles endowed with closed 7-forms. This notion was further extended to sphere bundles by Lind, Sati and Westerland \cite{LSW16} as a duality between $S^{2n-1}$-bundles endowed with closed $(4n-1)$-forms. We generalise this relation one step further and define T-duality for $S^{2n-1}$-bundles endowed with closed odd forms of arbitrary degree. The degree of the form determines the dimension of the fibers of the dual spaces. We show that $T$-duals exist and, as in the previous cases, $T$-dual spaces have isomorphic twisted cohomology. We finish by introducing a version of Courant algebroids which is compatible with spherical T-duality.}
\vskip12pt
\noindent
MSC classification 2020: 	55R25, 57R22, 57R19, 57N65.\\
Subject classification: Differential Topology.

\tableofcontents

\section{Introduction}

Target-space duality, or T-duality for short, is a duality that comes from physics \cite{Bus87,RV92}  in the presence of a torus symmetry but also has very concrete mathematical formulations and consequences. In its simplest form, when the target-space is a Riemannian circle, it is marked by inversion of the radius of the circle and swapping of physical quantities (winding and momentum) to yield equivalent physical theories.

Mathematically, T-duality consists of two ingredients. The first, {\it topological T-duality} \cite{BEM04}, relates the global topology of dual target-spaces: the background form on one side influences the topology of the dual space, there are isomorphisms between the twisted cohomologies of T-dual spaces and also of their twisted $K$-theories. The second ingredient, that makes T-duality {\it geometric}\footnote{Not standard nomenclature.}, is that topological duality is compatible with an isomorphism of Courant algebroids over T-dual spaces \cite{CG11,cavalcanti-2004}. This isomorphism is the mathematical tool that allows one to transport geometric structures between T-dual spaces, something that physicists had already been doing from the start, years before any of this was mathematically encoded.

The extension of T-duality to other setups has been a rich area of research in both physics and mathematics. From the perspective of transporting geometric structures and possible relations to physics, T-duality is a particular case of Poisson-Lie duality \cite{KlSe95}, can be extended to transitive Courant algebroids in this way mathematically encoding heterotic duality \cite{MR3418511,GARCIAFERNANDEZ20191059} and more \cite{Cortes:2021aab} and has been extended to a duality of Lie algebroids \cite{CaWi22}. From the topological perspective, common questions regard the extension to torus actions with fixed points \cite{MR2246781,CaWi22} and including non-Abelian symmetries \cite{DELAOSSA1993377}.

For the case of non-Abelian symmetries, Bouwknegt, Evslin and Mathai proposed in  \cite{BEM15} a version of {\it topological T-duality} which is quite distinct from other versions in the literature, e.g. Poisson-Lie T-duality. Here T-duality relates the total spaces of $\mathrm{SU}(2)$-bundles and, as in the original topological T-duality, switches information from a background closed form of fixed degree with the topology of the dual space. Once again, T-dual spaces have isomorphic twisted cohomologies. In this framework, the structure group, $\mathrm{SU}(2)$, is shared among T-duals instead of relating principal bundles with dual groups. A minor shortcoming is that even when the background form represents an integral class, there may be no T-dual space. This construction was further extended in \cite{BEM15b} to oriented 3-sphere bundles  and in \cite{LSW16} to $(2n-1)$-sphere bundles over $d$-dimensional manifolds. Up to now no corresponding version of {\it geometric T-duality} had emerged.

It is this formulation of \cite{BEM15,BEM15b,LSW16} that we will develop further. Here we allow for the background form to have any (fixed) odd degree higher than the dimension of the sphere fibers. Under T-duality, topological information of the background form is still exchanged with the Euler class of the dual bundle. With this framework we still have isomorphisms of twisted cohomology and one can construct a T-dual as long as the background form is (a multiple of) an integral class. The new feature is that while T-dual pairs fiber over a common base, the dimension of the fiber needs not to agree and that information is encoded in the degree of the background form.

In our study we single out a subspace of the space of forms that plays a similar role $S^1$-invariant forms did for usual T-duality. Focusing on this subspace, we propose an accompanying higher Courant algebroid for spherical T-duality and prove that T-dual spaces have isomorphic higher Courant algebroids, thus realising geometric spherical T-duality. Finally we relate the subbundle of self T-dual sections of this Courant algebroid with the underlying Courant algebroid that appears in $\mathrm{E}_6$ expectional generalized geometry.

This paper is organised as follows: in Section \ref{sec:t-duality} we give the definition of (topological spherical) T-duality and state its first basic properties. Section \ref{sec:angular form} recaps the definition and main properties of a global angular form and Euler class of a sphere bundle. With that form at hand, we can establish the existence theorem (Theorem \ref{theo:existence}) in Section \ref{sec:existence}. In Section \ref{sec:twisted} we recap the definiton of twisted cohomology and establish a key result of the twisted cohomology of sphere bundles (Theorem \ref{prop:quasi-isomorphism2}). In Section \ref{sec:isomorphism} we prove that T-dual spaces have isomorphic twisted cohomologies (Theorem \ref{theo:isomorphism}). In Section \ref{sec:examples} we give a few examples of T-dual spaces. In Section \ref{sec:geometric} we move to geometric T-duality, first  introducing the relevant version of Courant algebroids (Definition \ref{def:extended higher}) and then proving the T-duality isomorphism theorem (Theorem \ref{theo:geometric}).

\section{Spherical T-duality}\label{sec:t-duality}

Our object of interest are pairs $(E,H)$, where $E$ is an oriented sphere bundle, $\pi\colon E \to M$, with odd dimensional fibers and  $H\in \Omega^\bullet(E)$ is a closed odd-form of a specific degree. Throughout this paper we assume the base manifold to be connected. Given  a fibration $f\colon M^{n+d} \to N^n$ with compact oriented fibers, we denote by $f_*\colon \Omega^k(M) \to \Omega^{k-d}(N)$ integration along the fibers of $f$ and by $f^*\colon \Omega^k(N) \to \Omega^k(M)$ the pull-back of forms from $N$.

 Topological spherical T-duality, or just T-duality for short, is a relation between two such objects:

\begin{definition}
Let $n,k \in \NN\backslash\{0\}$, let $(E,H)$ be an oriented $(2n-1)$-sphere bundle over $M$ with a closed form $H \in \Omega^{2(n+k)-1}_{cl}(E)$ and let $(\hat{E},\hat{H})$ be an oriented $(2k-1)$-sphere bundle over $M$ with a closed form $\hat{H} \in \Omega^{2(n+k)-1}_{cl}(\hat{E})$. The bundles $(E,H)$ and $(\hat{E},\hat{H})$ are \emph{topological spherical T-duals} if there is a form $F\in \Omega^{2(n+k)-2}(E\times_M\hat{E})$ such that
\begin{itemize}
\item $dF = p^*H - \hat{p}^*\hat{H}$, where
$p\colon E\times_M\hat{E} \to E$ and $\hat{p}\colon E\times_M\hat{E}$ are the natural projections and
\item There is $x \in M$  for which $\pi_* \circ p_* (F): = \int_{E_x\times \hat{E}_x} F \neq 0$.
\end{itemize}
T-dual spaces $(E,H)$ and $(\hat{E},\hat{H})$ form a \emph{unimodular T-dual pair} if
 \[\pi_* \circ p_*(F)  =1.\]
\end{definition}

We gather the information on T-dual spaces in the following commutative diagram
\begin{equation}\label{eq:basic setup}
\xymatrix@C=-18pt{
& (E\times_M \hat{E}, dF = p^*H - \hat{p}^*\hat{H})\ar[dl]^{p}\ar[dr]_{\hat{p}}&\\
(E,H)\ar[dr]^{\pi}& & (\hat{E},\hat{H})\ar[dl]_{\hat{\pi}}\\
& M&
}
\end{equation}

We make a couple of simple remarks about T-duals.

\begin{lemma}Given a T-dual pair $(E,H)$ and $(\hat{E},\hat{H})$ over a manifold $M$, the function $\pi_* \circ p_*(F)\colon M \to \R$ is constant and nonvanishing.
\end{lemma}
\begin{proof}
Indeed, the integral under consideration is both the fiber integration $(\pi \circ p)_*$ and $-(\hat\pi \circ \hat{p})_*$ and since the degree of $F$ is the dimension of the fibers, $(\pi \circ p)_* F \in \Omega^0(M)$ and we have
\[d ((\pi \circ p)_* F) = (\pi \circ p)_* d F =  (\pi \circ p)_* (p^*H- \hat{p}^*\hat{H}) = (\pi \circ p)_* p^*H +  (\hat{\pi} \circ \hat{p})_* \hat{p}^*\hat{H} =0,\]
showing that the function $(\pi \circ p)_* F$ is constant and by assumption it does not vanish in at least one point.
\end{proof}

\begin{lemma}\label{lem:H vs H'}
 If $(E,H')$ and $(\hat{E},\hat{H}')$ are (unimodular) T-duals, $H'$ is cohomologous to $H$ and $\hat{H}'$ is cohomologous to $\hat{H}$, then $(E,H)$ and $(\hat{E},\hat{H})$ are also (unimodular) T-duals.
\end{lemma}
\begin{proof}
Indeed, if $p^*H' - \hat{p}^*\hat{H} = dF$ (because $(E,H')$ and $(\hat{E},\hat{H}')$ are T-duals), $H - H' = dB$ and $\hat{H}-\hat{H}' = d\hat{B}$, then $F + p^*B -\hat{p}^*\hat{B}$ expresses $(E,H)$ and $(\hat{E},\hat{H})$ as T-duals.
\end{proof}

Finally the relationship between the cohomology classes $p^*H$ and $\hat{p}^*\hat{H}$ in the correspondence space translates into specific relations on $M$ and $\hat{M}$:
\begin{lemma}\label{lem:euler class becomes exact}
Given a T-dual pair $\pi\colon (E,H) \to M$ and $\hat{\pi}\colon (\hat{E},\hat{H}) \to M$ we have
\[[\hat{\pi}^* \pi_* H] = 0 \in H^\bullet(\hat{M}) \quad \mbox{ and } \quad [\pi^*\hat{\pi}_* \hat{H}] = 0 \in H^\bullet(M).\]
\end{lemma}
\begin{proof}
Given that T-duality is a symmetric relation, it is enough to prove the first claim, which follows from a simple computation:
\[0 = [\hat{p}_*dF] =  [\hat{p}_*(p^*H - \hat{p}^*\hat{H})] = \hat{p}_*(p^*H)= [\hat{\pi}^* \pi_* H].\]
\end{proof}


\section{Global angular form}\label{sec:angular form}
The key ingredient for the construction of T-dual spaces, proof of the isomorphism of twisted cohomologies and eventual introduction of a geometric interpretation of T-duality is a global angular form, which is carefully explained, for example, in \cite{bott-tu}. We recall below its definition and main properties. 

\begin{definition}
Given an oriented sphere bundle $\pi\colon E \to M$ with fiber $S^m$, a {\it global angular form} on $E$ is a form $\psi \in \Omega^m(E)$ such that
\begin{itemize}
\item There is an $x \in M$,  for which $\pi_*\psi|_x:= \int_{E_x} \psi = 1$, where $E_x$ is the fiber over $x$,
\item $d\psi = \pi^*\varepsilon$, for some form $\varepsilon \in \Omega^{m+1}(M)$.   
\end{itemize}
\end{definition}
The second condition implies that $\pi_*\psi$ is a constant function and hence if the normalization of $\psi$ holds at one fiber it holds at all fibers.

\begin{remark}
One can further require that for all $x\in M$ the restriction of the global angular form to  $E_x$ is in fact a volume form. As it turns out, we will not need that assumption but imposing it is also not a restriction, since any oriented sphere bundle admits a global angular form which also satisfies this condition \cite{bott-tu}. It may quite well be that for specific geometric applications and connections to physics not explored in this paper, one must be careful about the allowed global angular forms.
\end{remark}
\begin{lemma}[Bott--Tu \cite{bott-tu}, pg 118-119]
Any two choices of global angular form on an oriented sphere bundle $E \to M$ give rise to the same cohomology class $\mc{E} = [\varepsilon]\in H^{m+1}(M)$.
\end{lemma}

\begin{definition}
The class $\mc{E}$ defined above is the \emph{Euler class} of the oriented sphere bundle.
\end{definition}

Given a sphere bundle $\pi \colon E \to M$, fiber integration, $\pi_*\colon \Omega^{k+m}(E) \to \Omega^k(M)$, is characterized by the property
\[\pi_*(\psi \wedge (\pi^*\gf)) = \gf,\]
where $\psi$ is any global angular form and $\gf \in \Omega^k(M)$.

Throughout this paper we will use a fundamental result about oriented sphere bundles that regards a subcomplex of the complex of differential forms:

\begin{definition}
Let $E \to M$ be an oriented sphere bundle and let $\psi$ be a global angular form. The \emph{Gysin complex} associated to $E$ and $\psi$ is the graded algebra
\begin{equation}\label{eq:subcomplex}
\Omega_\psi^{\bullet} (E) :=\{\pi^*\gf_0 + \psi \wedge \pi^*\gf_1 \colon \gf_i\in \Omega^\bullet(M)\} = \Gamma(M; \wedge\IP{\psi}\tensor \wedge^\bullet T^*M)
\end{equation}
\end{definition}

Not only is the Gysin complex is a subalgebra of $\Omega^\bullet(E)$, but it is also a subcomplex.

\begin{proposition}\label{prop:quasi-isomorphism}
The Gysin complex forms a subcomplex of $(\Omega^\bullet(E),d)$ and its inclusion $(\Omega_\psi^{\bullet}(E),d)  \hookrightarrow (\Omega^\bullet(E),d)$ is a quasi-isomorphism.
\end{proposition}
\begin{proof}
The fact that the subspace $\Omega_\psi^{\bullet} (E)$ is preserved by the exterior derivative follows from the fact that $d\psi$ is basic. The fact that the inclusion is a quasi-isomorphism is just the Gysin sequence re-dressed. Indeed, it is immediate that we have a short exact sequence of differential complexes
\[0 \to \Omega^\bullet(M) \stackrel{\pi^*}{\to}\Omega_\psi^{\bullet} (E) \stackrel{\pi_*}{\to} \Omega^{\bullet-m}(M) \to 0,\]
and one readily sees that the corresponding long exact sequence
\[\cdots \to H^\bullet(M) \stackrel{\pi^*}{\to}H(\Omega_\psi^{\bullet}(E)) \stackrel{\pi_*}{\to} H^{\bullet-n}(M) \stackrel{\cup \mc{E}}{\to} H^{\bullet+1}(M) \to \cdots,\]
is exactly the Gysin sequence that computes the de Rham cohomology of $E$ and the claim follows from the five-lemma.
\end{proof}

If $m$ is even, then $\mc{E} =0$, and this is one of the reasons why we will focus our attention on sphere bundles with odd dimensional fibers, say $m = 2n-1$. 

\section{Existence}\label{sec:existence}

With the global angular form and the Euler class at hand we are ready to establish the (in)existence of T-duals. This is done in a constructive way. The first step is to establish the obstruction to the existence of T-duals.
\begin{proposition}\label{prop:dual euler}
Given a T-dual pair $(E,H)$ and $(\hat{E},\hat{H})$ over a manifold $M$, the class $[\pi_* H] \in H^\bullet(M)$ is a multiple of the Euler class of $\hat{E}$. In particular, given an oriented sphere bundle with closed form $(E,H)\to M$ if the class $[\pi_* H] \in H^\bullet(M)$ is not a multiple of an integral class then $(E,H)$ has no T-dual.
\end{proposition}
\begin{proof}
The proposition holds trivially in the case  $[\pi_* H] = 0$, so we focus on the case $[\pi_* H] \neq 0$. It follows from the Gysin sequence for the cohomology of sphere bundles that the only way that one can have $[\hat{\pi}^*\pi_*H] = 0$ (as in Lemma \ref{lem:euler class becomes exact}), is if $[\pi_*H]$ is a multiple of the Euler class of $\hat{E}$. Since a cohomology class is a multiple of an Euler class if and only if it is a multiple of an integral class \cite{BV03}, the last claim follows.
\end{proof}

This integrality condition on $[\pi_*H]$ turns out to be the only obstruction to the existence of T-duals:

\begin{theorem}[Existence of T-duals]\label{theo:existence}
Let $\pi\colon E\to M$ be an oriented $S^{2n-1}$-bundle over a manifold $M$ and let $H \in \Omega_{cl}^{2(n+k)-1}(E)$ be a closed form with $n,k>0$. If the fiber integration $\pi_*[H] \in H^{2k}(M)$
is a multiple of an integral class, then there is a $(2k-1)$-sphere bundle $(\hat{E},\hat{H})$ T-dual to $(E,H)$. If further $\pi_*[H]$ is the Euler class of a sphere bundle over $\hat{E} \to M$, then  there is a sphere bundle $(\hat{E},\hat{H})$ for which $(E,H)$ and $(\hat{E},\hat{H})$ are a unimodular T-dual pair. 
\end{theorem}
\begin{proof}
Pick a global angular form $\psi$ for $E$ with $d\psi = \pi^*\varepsilon$ the corresponding representative for the Euler class of $E$.

Due to Lemma \ref{lem:H vs H'} and Proposition \ref{prop:quasi-isomorphism}, we may assume that $H \in \Omega^{2(n+k)-1}_{\psi}(E)$, say $H = \pi^*H_1 + \psi \wedge \pi^* H_0$. Then $H$ being closed is equivalent to
\begin{equation}\label{eq:dh=0}
\begin{aligned}
dH_0 &= 0\\
dH_1 + \varepsilon  \wedge  H_0 &=0
\end{aligned}
\end{equation}

Since $\pi_*H = H_0$ represents a multiple of an integral class, there is $\lambda \in \R_+$ such that $\lambda [H_0]$ is the Euler class, $\hat{\mc{E}}$, of a $(2k-1)$-sphere bundle $\hat{\pi}\colon \hat{E}\to M$ (see for example \cite[Lemma B1]{BV03}). Notice that if $H_0$ is already the Euler class of a sphere bundle we may use $\lambda =1$.

Let $\hat{\psi}$ be a global angular form for $\hat{E}$. After changing $\hat{\psi}$ by the pullback of an exact form we can arrange that $d\hat{\psi} = \lambda \hat\pi^*H_0 =: \hat\pi^*\hat{\varepsilon}$. We endow $\hat{E}$ with the form
\[\hat{H} = \hat{\pi}^*H_1 + \frac{1}{\lambda}\hat{\psi} \wedge \hat{\pi}^*\varepsilon.\]
For these choices, using \eqref{eq:dh=0}, we have
\[d\hat{H} = \hat{\pi}^*dH_1 + \frac{1}{\lambda}d\hat{\psi} \wedge \hat{\pi}^*\varepsilon = \hat{\pi}^*(dH_1 + H_0 \wedge \varepsilon) = 0,\]
showing that $\hat{H}$ is a closed form. Further, on $E \times_M \hat{E}$ we have
\begin{align*}
p^*H - \hat{p}^*\hat{H}& =  p^*\pi^* H_1 + \frac{1}{\lambda}p^*\psi \wedge p^* \pi^* \hat{\varepsilon} - \hat{p}^*\hat{\pi}^*H_1 - \frac{1}{\lambda}\hat{p}^*\hat{\psi} \wedge \hat{p}^*\hat{\pi}^*\varepsilon\\
&=  \frac{1}{\lambda}p^*\psi \wedge \hat{p}^*d\hat{\psi} -  \frac{1}{\lambda}\hat{p}^*\hat{\psi} \wedge \hat{p}^*d\psi\\
&=  \frac{1}{\lambda}d(\hat{p}^*\hat{\psi}\wedge p^* \psi).
\end{align*}
Hence we can take $F = \frac{1}{\lambda}\hat{p}^*\hat{\psi}\wedge p^* \psi$, which satisfies the normalization
\[\pi_* \circ p_*F = \frac{1}{\lambda} \neq 0,\]
showing that $(E,H)$ and $(\hat{E},\hat{H})$ are a T-dual pair and that if $\hat{E}$ itself was already the Euler class of a sphere bundle, $\hat{E}$, then we could have taken $\lambda =1$ and obtain a unimodular pair.
\end{proof}

Notice that there were several choices in the construction above: Firstly, the choice a multiple of $\hat{\mc{E}}$ which is the Euler class of a sphere bundle. Secondly the choice of a sphere bundle with that specific Euler class. Here not only the presence of torsion in integral cohomology brings ambiguity but also that fact that in general the integral Euler class does not determine the sphere bundle. All these choices affect the diffeomorphism type of $\hat{E}$. Imposing unimodularity only adresses some of these issues. So, in general, T-duals are not unique, a fact that was also explored in \cite{BEM15b} to produce infinitely many non diffeomorphic T-duals. 

%
%

\section{Twisted cohomology}\label{sec:twisted}

Our next objective is to show that there is a clear topological relation between T-dual spaces. Any such relation must take the background form into consideration. In our case we consider the twisted cohomology by the form $H$.
\begin{definition}
The \emph{twisted cohomology} of a manifold endowed with a closed odd-form, $(E,H)$, is the cohomology of the operator $d^H = d + H \wedge$.
\end{definition}
 There are two points worth of notice:
\begin{itemize}
\item For a form $H\in \Omega^{2m+1}(E)$, the $d^H$ is only compatible with degree modulo $2m$, therefore the $d^H$-cohomology is only $\Z_{2m}$-graded.
\item The isomorphism type of the $d^H$-cohomology (and many other objects one can construct depending on $H$) depends only on the cohomology class of the form $H$:
\end{itemize}
\begin{lemma}\label{lem:lemma5}
Given two cohomologous forms $H$, $H' \in \Omega^{2m+1}(E)$, let $B$ be such that $H - H'= dB$. Then the following is an isomorphism of differential complexes
\[(\Omega^\bullet(E),d^H) \stackrel{e^B}{\longrightarrow}  (\Omega^\bullet(E),d^{H'}),\]
in particular these complexes have isomorphic cohomologies.
\end{lemma}

We can put Proposition \ref{prop:quasi-isomorphism} and Lemma \ref{lem:lemma5} together to find a subspace of $\Omega^\bullet(M)$ whose cohomology is quasi-isomorphic to the $d^H$-cohomology for sphere bundles.

\begin{proposition}\label{prop:quasi-isomorphism2}
Let $(E,H') \to M$ be an oriented $(2n-1)$-sphere bundle with global angular form $\psi$ and $H' \in \Omega^{2m+1}(E)$. Let $H\in \Omega_\psi^{\bullet} (E)$ be a closed form cohomologous to $H'$, say $H-H' =dB$. Then $\Omega_\psi^{\bullet} (E)$ is preserved by $d^{H}$ and the following maps are quasi-isomorphisms
\[(\Omega_\psi^{\bullet} (E),d^{H})\stackrel{\iota}{\hookrightarrow} (\Omega^{\bullet} (E),d^{H}) \stackrel{e^{B}}{\longrightarrow}(\Omega^{\bullet} (E),d^{H'}).\]
\end{proposition}
\begin{proof}
The second map is a isomorphism by Lemma \ref{lem:lemma5}, so we only need to prove that $(\Omega_\psi^{\bullet} (E),d^{H})$ is a differential complex and the inclusion $\iota\colon(\Omega_\psi^{\bullet} (E),d^{H})\hookrightarrow (\Omega^{\bullet} (E),d^{H})$ is a quasi-isomorphism. 

Since the fibers of $E$ are odd-dimensional, $\psi^2 =0$ and $\Omega_\psi^{\bullet} (E)$ is in fact a subalgebra and hence $d^{H}$ also preserves $\Omega_\psi^{\bullet} (E)$ making it a subcomplex of $(\Omega^{\bullet} (E),d^H)$. The general argument to show that the inclusion is a quasi-isomorphism is based on the fact that the $d^H$-cohomology of both $\Omega_\psi^{\bullet} (E)$ and $\Omega^{\bullet} (E)$ can be computed via a spectral sequence (see, for example, \cite[Section 1.4]{cavalcanti-2004}) whose first page is the de Rham cohomology of the corresponding complex. Since, by Proposition \ref{prop:quasi-isomorphism}, the inclusion is a quasi-isomorphism of de Rham cohomologies, the inclusion gives an isomorphism between the first page of these spectral sequences and therefore it also gives an isomorphism of their last pages. Below we spell part of this argument out following a tic-tac-toe approach.

To show that the inclusion is a quasi-isomorphism we need to show that the corresponding map in cohomology is injective and surjective. We will only prove surjectivity as injectivity is proved in a very similar way.

Given a cohomology class in $H^\bullet(\Omega^{\bullet} (E),d^H)$ we pick a representative $\gf$ for it and show that $\gf$ is cohomologous to a form in the image of the inclusion. We construct the form in the image of $\iota$ inductively on the degree of the homogenous components of $\gf$, starting from the lowest degree component and moving inductive up in degree. We split $\gf = \gf_0 + \gf_1 +\dots + \gf_l$, with $\gf_0$ a form of degree $p$ and $\gf_i$ a form of degree $p + 2m i$.

The condition $d^H \gf = 0$ translates into a series of conditions:
\[ d\gf_0 = 0, \qquad  d\gf_{i} + H\wedge \gf_{i-1} =0.\]
The first step on the induction is to change $\gf_0$ into a form in the image of the inclusion. Since $\gf_0$ is closed, there is a $\xi_0 \in \Omega^p_{\psi}(E)$ is the same cohomology class of $\gf_0$, that is $\gf_0 = \xi_0+ d\eta_0$. Consider $\gf - d^H\eta_0$. This form still represents the same $d^H$-cohomology class as $\gf$ but its lower degree term lies in  $\Omega^\bullet_{\psi}(E)$.

Assuming inductively that $\gf_i \in \Omega^\bullet_{\psi}(E)$ for $i <i_0$, we show that we can find another form cohomologous to $\gf$ for which $\gf_i \in \Omega^\bullet_{\psi}(E)$ for $i  < i_0+1$. We consider the equation
\[H\wedge \gf_{i_0-1} = - d\gf_{i_0}\]
Since $H$ and $\gf_{i_0-1} \in \Omega^\bullet_{\psi}(E)$ their product is also in that space and by the equation above it represents the zero cohomology class. By Proposition \ref{prop:quasi-isomorphism}, there is an element $\xi_{i_0}\in   \Omega^\bullet_{\psi}(E)$ such that
\[H\wedge \gf_{i_0-1} = - d\xi_{i_0}.\]
Therefore $d(\gf_{i_0} - \xi_{i_0}) =0$ and there is a form $\xi_{i_0}' \in  \Omega^\bullet_{\psi}(E)$ in the same cohomology class, so
\[ \gf_{i_0} - \xi_{i_0} - \xi_{i_0}'  = d\eta_{i_0},\]
for some $\eta_{i_0} \in \Omega^\bullet(E)$.

Consider now $\gf - d^H\eta_{i_0}$. This form is in the same $d^H$-cohomology class as $\gf$, the  terms of degree up to  
$2m(i_0-1)+p$ agree with those of $\gf$ and in degree $2mi_0+p$ we have
\[\gf_{i_0} - d\eta_{i_0} =  \xi_{i_0} +\xi_{i_0}' \in \Omega_{\psi}^\bullet(E),\]
thus completing the induction step.
\end{proof}

\section{The isomorphism of twisted cohomologies}\label{sec:isomorphism}
Now that we know that T-dual spaces exist and it is actually easy to construct T-dual pairs, we move on to the question of what the relation between these spaces is. Here we show one such relation:

\begin{theorem}\label{theo:isomorphism} 
Let $(E,H)$ and $(\hat{E},\hat{H})$ be a T-dual pair with $dF = p^*H - \hat{p}^*\hat{H}$. Then the map
\begin{equation}\label{eq:tau}
\tau_F \colon (\Omega^\bullet(E),d^H) \to (\Omega^\bullet(\hat{E}),d^{\hat{H}}),\qquad \tau_F(\gf) = \hat{p}_*e^Fp^*\gf.
\end{equation}
induces an isomorphism of twisted cohomologies which shifts degrees up by $2k-1$ modulo $2(n+k-1)$.
\end{theorem}

\begin{remark}
The map that realises the isomorphism in cohomology is neither injective nor surjective at the cochain level. 
\end{remark}
\begin{proof}[Start of the proof of Theorem \ref{theo:isomorphism}]
For any $F \in \Omega^{2(n+k-1)}(E\times_M\hat{E})$ we can form the map \eqref{eq:tau}. If $dF = p^*H -\hat{p}^*\hat{H}$, this is a map of cochain complexes, since
\[\tau_{F}(d^H\gf) = \hat{p}_*e^F p^* d^{H}\gf = \hat{p}_*d^{p^*H- dF}(e^F p^*\gf) =   \hat{p}_*d^{\hat{p^*}\hat{H}}(e^F p^*\gf) = d^{\hat{H}}\hat{p}_*e^F p^*\gf = d^{\hat{H}}\tau_F(\gf).\]

Our task is to show that under the nondegeneracy condition, $(\pi\circ p)_*F\neq 0$, this map is a quasi-isomorphism. We will achieve this by relating $\tau_F$ to a collection of other maps all of which are quasi-isomorphisms if and only if $\tau_F$ itself is. Eventually we will be able to compare $\tau_F$ to a map that is simple enough that we can check it is a quasi-isomorphism directly. Part of the difficulty is that $F$ at the moment is too general a form for us to get our hands on $\tau_F$ concretely.
\noqed
\end{proof}

\begin{lemma}\label{lem:exact change} In the situation above, if $F' = F + dB$ then $\tau_{F'} =  \hat{p}_*e^{F'}p^*$ induces the same map as $\tau_F$ in twisted cohomology.
\end{lemma}
\begin{proof}
We observe that
\[\tau_{F'} (\gf) - \tau_F(\gf) =  \hat{p}_*(e^{F'}-e^{F})p^*\gf =  \hat{p}_*(e^{dB}-1)e^{F}p^*\gf.\]
If $d^H \gf =0$, then $d^{\hat{p}^*\hat{H}}(e^F\gf)=0$ and hence
\[(e^{dB}-1)e^{F}p^*\gf = (dB)(\sum_{k=1}^\infty \frac{(dB)^{k-1}}{k!})e^{F}p^*\gf = d^{\hat{p}^*\hat{H}}\left(B (\sum_{k=1}^\infty \frac{(dB)^{k-1}}{k!})e^{F}p^*\gf\right),\]
showing that the difference $\tau_{F'} (\gf) - \tau_F(\gf)$ is $d^{\hat{H}}$-exact and hence $\tau_{F'}$ and $\tau_{F}$ induce the same map in twisted cohomology.
\end{proof}

\begin{lemma}\label{lem:lemma2}
Let $(E,H)$ and $(\hat{E},\hat{H})$ be a T-dual pair. If $F$ and $F'$ satisfy $dF = dF' = p^*H -\hat{p}^*\hat{H}$ and $(\pi\circ p)_*F' = (\pi\circ p)_*F\neq 0$, then there are closed forms $F_1 \in \Omega^{2(n-k-1)}(E)$  and $F_2 \in \Omega^{2(n-k-1)}(\hat{E})$ such that
\[\tau_{F} = e^{-F_2} \circ  \tau_{F'}\circ e^{F_1}.\]
In particular, $\tau_{F'}$ induces an isomorphism of twisted cohomologies if and only if  $\tau_F$ does.
\end{lemma}
\begin{proof}
Under the hypothesis, $F - F'$ is a closed $2(n+k-1)$-form for which $(\pi\circ p)_*(F - F')= 0$. Using the spectral sequence for the cohomology of the double sphere bundle, $E\times_M \hat{E}\to M$,  we have that $F - F'$ is cohomologous to $p^*F_1 - \hat{p}^*F_2$ for some closed forms $F_1 \in \Omega^{2(n-k-1)}(E)$ and  $F_2 \in \Omega^{2(n-k-1)}(\hat{E})$. In particular $F- F' - p^*F_1 +\hat{p}^*F_2$ is exact, hence by Lemma \ref{lem:exact change} and $\tau_F$ and $\tau_{F'+ p^*F_1  - \hat{p}^*F_2}= e^{-F_2} \circ  \tau_{F}\circ e^{F_1} $ induce the same map in twisted cohomology.

Since the maps $e^{-F_2}$ and $e^{F_1}$ are isomorphisms of twisted cohomology, we conclude that $\tau_{F'}$ is an isomorphism of twisted cohomology if and only if $\tau_F$ is.
\end{proof}

Now we start to move towards making a concrete simplification to the map $\tau_F$.
\begin{lemma}\label{lem:F in Gysin complex}
Let $(E,H)$, $(\hat{E},\hat{H})$ be a T-dual pair with $F = p^*H - \hat{p}^*\hat{H}$. Let  $\psi$ and $\hat{\psi}$ be global angular forms for $E$ and $\hat{E}$. There is a form $F'$ in the Gysin complex of $E\times_M \hat{E}$
\[\Omega_{\psi,\hat{\psi}}(E\times_M \hat{E}) = \{(\pi \circ p)^*\gf_0 + p^*(\psi \wedge \pi^*\gf_1) +  \hat{p}^*(\hat{\psi}\wedge \hat{\pi}^*\gf_2) + \hat{p}^*\hat{\psi} \wedge p^*\psi \wedge (\pi \circ p)^*\gf_3\colon \gf_i \in \Omega(M)\}\]
such that $\tau_F$ an $\tau_{F'}$ induce the same map in twisted cohomology.
\end{lemma}
\begin{proof}
Due to Proposition \ref{prop:quasi-isomorphism} and Lemma \ref{lem:H vs H'}, we may assume that $H \in \Omega_\psi^\bullet(E)$ and $\hat{H} \in \Omega_\psi^\bullet(\hat{E})$.  Due to Proposition \ref{prop:quasi-isomorphism} the inclusion $\Omega_{\psi,\hat{\psi}}(E\times_M \hat{E}) \hookrightarrow\Omega(E\times_M \hat{E})$ is quasi-isomorphism.

Since on $E\times_M \hat{E}$, $H - \hat{H}$ is exact, there is a form $\tilde{F} \in \Omega^{2(n+k-1)}_{\psi,\hat{\psi}}(E\times_M \hat{E})$,  such that $d\tilde{F} = p^*H - \hat{p}^*\hat{H}$. Then $F - \tilde{F}$ is a closed form on $E\times_M \hat{E}$ and hence $(\pi\circ p)_*(F - \tilde{F})$ is constant on $M$. There are two cases to consider:

\vskip6pt
\noindent
{\it Case 1: $(\pi\circ p)_*(F - \tilde{F}) \neq 0$}. In this case, $F - \tilde{F}$ is a global closed form on $E \times_M \hat{E}$ that restricts fiberwise to a nonvanishing top degree cohomology class. From Gysin sequences for the cohomology of the $S^{2n-1}\times S^{2k-1}$-bundle, we see that this can happen if and only if the Euler classes of both $E$ and $\hat{E}$ vanish which is equivalent to the existence of closed global angular forms, $\psi$ and $\hat{\psi}$. In this case, if $(\pi\circ p)_*(F - F') =\lambda$ consider instead $\tilde{\tilde{F}}= \tilde{F}-\lambda \psi \wedge \hat\psi$. This form also satisfies $d\tilde{\tilde{F}} = p^*H - \hat{p}^*\hat{H}$ and $(\pi\circ p)_*(F- \tilde{\tilde{F}}) =0$, which puts us in the second case.

\vskip6pt
\noindent
{\it Case 2: $(\pi\circ p)_*(F - \tilde{F}) =0$}. From Lemma \ref{lem:lemma2}, we have that
\[\tau_{F} = e^{-F_2}\circ  \tau_{\tilde{F}} \circ e^{F_1},\]
for closed forms $F_1$ and $F_2$. By Lemma \ref{lem:lemma5}, changing $F_1$ and $F_2$ by exact elements does not change the map induced in cohomology and by Proposition \ref{prop:quasi-isomorphism}, we can choose $F_1$ and $F_2$ in the corresponding  Gysin complexes. Then the form $F' = \tilde{F} +p^*F_1-\hat{p}^*F_2$ satisfies the conditions of the lemma.
\end{proof}

\begin{lemma}\label{lem:T-duality of subcomplexes}
The map $\tau_{F'}$, with $F' \in\Omega_{\psi,\hat{\psi}}(E\times_M \hat{E})$, filters through $\Omega_\psi(E)$ and $\Omega_{\hat{\psi}}(\hat{E})$, that is, there is a map $\tau'_{F'} \colon \Omega_\psi(E) \to \Omega_{\hat{\psi}}(\hat{E})$ for which the following diagram commutes
\begin{equation}\label{eq:the diagram}
\begin{gathered}
\xymatrix{
(\Omega^\bullet(E),d^H) \ar[r]^{\tau_{F'}} & (\Omega^\bullet(\hat{E}), d^{\hat{H}})\\ 
(\Omega^\bullet_\psi(E),d^H) \ar[r]^{\tau'_{F'}}\ar[u]^\iota & (\Omega^\bullet_{\hat{\psi}}(\hat{E}), d^{\hat{H}})\ar[u]^{\hat{\iota}}
}
\end{gathered}
\end{equation}
Further $\tau'_{F'}$ is an isomorphism of differential complexes.
\end{lemma}
The fact that we obtain a map $\tau'_{F'}$ with the stated properties is a direct computation.

\begin{proof}[End of proof of Theorem \ref{theo:isomorphism}]
We let
$F' =  \hat{p}^*(\hat{\psi}\wedge \pi^*(F_1)) + \lambda p^*\psi \wedge \hat{p}^*\hat{\psi} + p^* (\psi \wedge \pi ^* (F_2)) + p^*\pi^*(F_3)$ with $F_i \in \Omega^\bullet(M)$ and $\lambda  = \hat{\pi^*}\circ \hat{p}^* F\neq 0$. Then we have
\begin{align*}
\tau_{F'}(\gf) &= \hat{p}_*\left( e^{\hat{p}^*(\hat{\psi}\wedge \pi^*(F_1)) + \lambda p^*\psi \wedge \hat{p}^*\hat{\psi} + p^* (\psi \wedge \pi ^* (F_2)) + p^*\pi^*(F_3)}  p^*\gf\right)\\
& =  e^{\hat{\psi}\wedge \pi^*(F_1)} \wedge \left(\hat{p}_*e^{\lambda p^*\psi \wedge \hat{p}^*\hat{\psi}}p^*\right)\left(e^{\psi \wedge \pi ^* (F_2)+ \pi^*(F_3)} \gf)\right),
\end{align*}
which expresses $\tau_{F'}$ as the composition of three maps:
\begin{itemize}
\item On the left, the automorphism of $\Omega_{\hat{\psi}}(\hat{E})$,
\[e^{\hat{\psi}\wedge \pi^*(F_1)}\colon \Omega_{\hat{\psi}}(\hat{E}) \to \Omega_{\hat{\psi}}(\hat{E}),\]
\item On the right, the automorphism  of $\Omega_{\psi}(E)$,
\[e^{\psi\wedge \pi^*(F_2) +\pi^*(F_3)}\colon \Omega_{\psi}(E) \to \Omega_{\psi}(E),\]
\item And in the middle the map
\begin{equation}
\tau'_{\lambda p^*\psi \wedge \hat{p}^*\hat{\psi}}(\gf) :=  \hat{p}_*e^{\lambda p^*\psi \wedge \hat{p}^*\hat{\psi}} p^* \gf
\end{equation}
\end{itemize}
Since the outer maps are automorphisms, to prove the existence of the map $\tau'_{F'}$ (and that it is an isomorphism) we only need to show that the middle map is an isomorphism. That we can do very explicitly for $\gf = \pi^*\gf_0 + \psi \wedge \pi^*\gf_1$:
\begin{align*}
\tau'_{\lambda p^*\psi \wedge \hat{p}^*\hat{\psi}}(\pi^*\gf_0 + \psi \wedge \pi^*\gf_1) &=  \hat{p}_*e^{\lambda p^*\psi \wedge \hat{p}^*\hat{\psi}} p^*(\pi^*\gf_0 + \psi \wedge \pi^*\gf_1)\\
& = \hat{p}_*(p^*\pi^*\gf_0 + p^*\psi \wedge p^* \pi^*\gf_1 + \lambda p^*\psi \wedge \hat{p}^*\hat{\psi}\wedge p^*(\pi^*\gf_0))\\
& = \hat{p}_*(p^*\pi^*\gf_0 + p^*\psi \wedge \hat{p}^* \hat{\pi}^*\gf_1 + \lambda p^*\psi \wedge \hat{p}^*\hat{\psi}\hat{p}^*(\hat{\pi}^*\gf_0))\\
& = \hat{\pi}^*\gf_1 + \lambda  \hat{\psi}\wedge \hat{\pi}^*\gf_0,
\end{align*}
showing that $\tau'_{\lambda p^*\psi \wedge \hat{p}^*\hat{\psi}}$ is an isomorphism and hence so is $\tau'_{F'}$.

Finally the claim that $\tau_{F'}$ itself is a quasi-isomorphism follows from the fact that the remaining three maps in diagram \eqref{eq:the diagram} are.
\end{proof}

\subsection{$K$-theoretical interpretation}

Another common feature of topologically T-dual spaces is that they have isomorphic twisted K-theories \cite{BEM15,LSW16}. This still holds for the present case, with the same proof as in \cite{LSW16}. We finish this section with an indication of the steps needed to arrive at such isomorphism. A vague statement of the result is:

\begin{theorem}
Let $(E,H) \to M$ and $(\hat{E},\hat{H})\to M$ be a unimodular T-dual pair for which $H$ and $\hat{H}$ are induced by $(2(n+k)-3)$-gerbes on $E$ and $\hat{E}$, respectively. Then the Fourier--Mukai transform on $E\times_M \hat{E}$ gives rise to an isomorphism of twisted K-theories.
\end{theorem}

Making sense of the notion of higher gerbes and higher twisted K-theory requires some effort and we refer to \cite{LSW16} for a precise description of the concepts involved. The point to notice is that the proof of the result in that paper relies on three properties of T-duality which still hold in the present context:
\begin{itemize}
\item (normalization) The result holds if the base manifold is a point (because of the unimodularity condition),
\item (naturality) T-duality is natural under pull-backs, that is, $(E,H)\to M$ and $(\hat{E},\hat{H})\to M$ are unimodular $T$-duals and $f\colon M'\to M$ is smooth then $(f^*E,f^*H)\to M'$ and $(f^*\hat{E},f^*\hat{H})\to M'$ are unimodular $T$-duals, and
\item (patching) T-duality preserves Mayer-Vietoris sequences.
\end{itemize}
Hence their proof still holds for the case in which the twisting form has degree $2(n+k-1)$.

Similarly, the higher Chern character still gives a natural transformation of twisted cohomology theories, between the twisted $K$-theory and twisted cohomology if one only considers the latter as $\Z_2$-graded \cite{LSW16}.

\section{Examples}\label{sec:examples}

\begin{example}
Consider the trivial sphere bundles $M \times S^{2n-1}$ and $M\times S^{2k-1}$ both endowed with the zero $(2(n+k)-1)$-form. Following the construction from Theorem \ref{theo:existence} these spaces are (unimodularly) T-dual to each other and hence, by Theorem \ref{theo:isomorphism}, their de Rham cohomologies (indexed modulo $2(n+k-1)$) are isomorphic via the map
\[a + [\sigma_{S^{2n-1}}] \cup b \mapsto  b + [\sigma_{S^{2k-1}}] \cup a, \qquad \mbox{ for }a,b \in H^\bullet(M),\]
where $\sigma_{S^{m}}$ denotes unit volume form on the sphere $S^m$. 
\end{example}

\begin{example}
Let $E \to M$ be a $(2n-1)$-sphere bundle over $M^{2n}$ with Euler class $\mc{E}$ and endowed with the zero form in $\Omega^{2(n+k)-1}(E)$. Then, following the construction of Theorem \ref{theo:existence}, $(E,0)$ is T-dual to $(M \times S^{2k-1}, \hat{\psi}\wedge \mc{E})$ for every $k$. That is, just as in the previous example, for the zero form there is ambiguity even on the dimension of the T-dual space.

Further, recall that there is a spectral sequence that computes the $d^{H}$-cohomology of $(E,H)$. The first page of this sequence is the de Rham cohomology with differential
\[[H]\cup \colon H^{\bullet}(E) \to  H^{\bullet}(E).\]
The second page of the spectral sequence that computes the $d^{H}$-cohomology is the cohomology of this this operator which we refer to as the $[H]$-cohomology. The following pages of the spectral sequence involve Massey products of the form
\[\IP{H,\dots,H, a}.\]
for $a \in \Omega_{cl}^\bullet(E)$.

In the present case, since on the $E$-side we took $H=0$, we see that the spectral sequence already degenerates at page $E_1$, while on the dual side, since $H \neq 0$, the page $E_2$ is different from $E_1$, showing that ``the index $n$ for which $E_n = E_\infty$ '' is not preserved by T-duality.
\end{example}

\begin{example}
Let $S^{2n+1} \to \C P^n$ be the Hopf fibration and endow $S^{2n+1}$ with $H= \sigma_{S^{2n+1}}$. Concretely, letting $\omega$ be the Fubini--Study form on $\C P^n$, we can endow $S^{2n+1}$ with a global angular form $\psi$ for which $d\psi = \pi^*\omega$ and $H = \psi \wedge \pi^* \omega^{n}$. Following the construction of Theorem \ref{theo:existence}, a T-dual to $(S^{2n+1}, \sigma_{S^{2n+1}})$ is $\hat{E}$, the $(2n-1)$-sphere bundle over $\C P^n$ with angular form $\hat{\psi}$ satisfying $d\hat{\psi} =\hat{\pi}^* \omega^n$ and background form $\hat{H} = \hat{\psi} \wedge \hat{\pi}^*\omega$.

The twisted cohomology of $(S^{2n+1},\sigma_{S^{2n+1}})$ vanishes hence so does the twisted cohomology of $(\hat{E},\hat{\psi} \wedge \hat{\pi}^*\omega)$, a fact that can also be checked directly, but with a few more steps.


\end{example}

\begin{example}
A more general version of the previous example, one can take  $(M^{2n},\omega)$ a compact K\"ahler manifold with K\"ahler form $\omega$ representing an integral class, and for sphere bundle we consider $\pi\colon E\to M$ the  principal circle bundle with Euler class $[\omega]$ and $H = \psi \wedge \pi^*\omega^n$.

For these choices, following the construction of Theorem \ref{theo:existence}, a T-dual for $E$ is the $(2n-1)$-sphere bundle $\hat{E}$ whose Euler class is $\omega^n$ endowed with the closed form $\hat{\psi} \wedge \hat{\pi}^*\omega$.

Since $H$ is a top-degree form on $E$, denoting by $[i]$ the reduction of $i$ modulo $2n$, the twisted cohomology of $E$ is
\[
H^{[i]}_{d^H}(E) = \begin{cases}H^i(E)&\qquad \mbox{ for } 0<i <2n,\\
0 & \qquad \mbox{ for } i = 0.
\end{cases}
\]
We can go one step further and compute the cohomology of $E$ in terms of the cohomology of $M$. Indeed, the Gysin sequence for $E\to M$ together with the Lefschetz property give
\[0 \to H^i(M) \stackrel{\cup[\omega]}{\hookrightarrow} H^{i+2}(M)\stackrel{\pi^*}{\to} H^{i+2}(E) \to 0 \qquad \mbox{ for } i< n-1,\]
\[0 \to H^i(E) \stackrel{\pi_*}{\to} H^{i-1}(M) \stackrel{\cup[\omega]}{\twoheadrightarrow}H^{i+1}(M) \to \{0\} \qquad \mbox{ for } i \geq n.\]
So the twisted cohomology of $(E,H)$ (and consequently also of $(\hat{E},\hat{H})$) is the primitive cohomology of $M$ up to middle dimension and its dual from middle dimension up.
\end{example}

\begin{example}[Self T-duality]
One situation that may be considered particularly interesting or intriguing is that of spaces which are self T-dual. This is because even though the space itself does not change under T-duality, the underlying T-duality map is potentially  a nontrivial isomorphism of the twisted cohomology of the space. In this example we spell out the conditions needed for a space to be self T-dual.

Of course, the first condition is that the dimension of the space and its T-dual should be the same. If $E \to M$ is a $2n-1$-sphere bundle, then this forces the background form $H$ to have degree $4n-1$. From Proposition \ref{prop:dual euler}, we have that $[\pi_* H]$ must be a multiple of the Euler class, $\mc{E}$ of $E$. The condition $dH =0$ implies that $\mc{E} \cup \mc{E} =0$.

Conversely, if 
\begin{align}
[\pi_* H] &= \lambda \mc{E} \label{eq:self T-dual 1}\\
\mc{E} \cup \mc{E} &=0. \label{eq:self T-dual 2}
\end{align}
Theorem \ref{theo:existence} shows that $(E,H)$ is self T-dual.

Even though T-duality is a property of the pair $(E,H)$, given a sphere bundle $S^{2n-1}\cdots E\to M$, one might be interested in knowing if there is an $H \in \Omega^{4n-1}(E)$ which makes $E$ self T-dual. In this case, condition \eqref{eq:self T-dual 2} is the only obstruction. Indeed, if  $\mc{E} \cup \mc{E} =0$, then, from the Gysin sequence,
\[\cdots H^{4n-1}(E) \stackrel{\pi_*}{\longrightarrow} H^{2n}(M) \stackrel{\cup \mc{E}}{\longrightarrow}H^{4n}(M) \to \cdots\]
we see that $\mc{E}$ arises as the fiber integration of a cohomology class $[H] \in H^{4n-1}(E)$. Taking $H$ any representative of this class makes $(E,H)$ unimodular self T-dual due to Theorem \ref{theo:existence}.

Notice that if the dimension of the base manifold is at most $4n-1$ (where the dimension of the sphere fiber is $2n-1$), then the condition $\mc{E} \cup \mc{E}$ holds trivially for dimensional reasons. In particular, if $\dim(E) \geq \frac{3}{2}\dim(M)$, there is a choice of background form on $E$ which makes it unimodular self T-dual.  
\end{example}

\section{Geometric Spherical T-duality}\label{sec:geometric}

As we saw in the previous sections, the question of existence of spherical T-duals and the immediate relationship between T-dual spaces is topological in nature. It is our biased point of view that for this new type of duality to have applications in differential geometry and physics it should be paired with an isomorphism of appropriate vector bundles.

In the original case of T-duality, for principal $S^1$-bundles, it was the introduction in \cite{cavalcanti-2004,CG11}  of a T-duality map, which is an isomorphism of {\it Courant algebroids over the base manifold} compatible with the cochain map \eqref{eq:tau}, that did the step of turning a topological relation into a geometric one. That step allowed for geometric structures to be transported between T-dual spaces and, for example, recovered the Buscher rules from physics in a differential geometric framework.

This passage from topological to geometric T-duality was missing for spherical T-duality up to this point. 
Since in the present case we are considering twists by closed forms of degree $2(n+k)-1$, a natural candidate for vector bundle for which there is an isomorphism would be the higher Courant algebroid $TE \oplus \wedge^{2(n+k)-3}T^*E$. But a direct try, even in the case of principal $\mathrm{SU}(2)$-bundles, shows that there is no T-duality isomorphism on the space of invariant sections of these bundles compatible with the map \eqref{eq:tau}. For general sphere bundles, where the notion of invariant form is absent, the path becomes even less clear.

We tackle this problem by first introducing a candidate of ``higher Courant algebroid'' that fits well with the current setup and then prove that T-duality indeed gives rise to an isomorphism of such Courant algebroids.

\subsection{Extended (higher) Courant algebroids}

With the introduction of a global angular form, the  Gysin complex \eqref{eq:subcomplex} and an explicit realisation of the quasi-isomorphism of differential complexes as an isomorphism of complexes (Lemma \ref{lem:T-duality of subcomplexes}), it becomes plausible that the Gysin complex \eqref{eq:subcomplex} is the relevant object for a geometric realisation. From this point of view we treat the global angular form, $\psi$, as a formal variable on $M$ extending the algebra of differential forms:
\[\wedge^\bullet_\psi T^*M := \wedge^\bullet\IP{\psi} \tensor \wedge^\bullet T^*M, \quad \mathrm{deg}(\psi) = 2n-1,\quad d\psi = \varepsilon.\]

In this context, a natural candidate for Courant algebroid compatible with T-duality would be an extension of $TM\oplus \wedge^{2(n+k)-3} T^*M$ obtained by adding the formal variable $\psi$ as above to the ``$T^*M$'' component and the dual variable, $\del_\psi$, which should correspond geometrically to a fiberwise top degree multivector field:
\[\mathrm{deg}(\del_\psi) = -2n+1,\quad  \IP{\psi,\del_\psi} = 1, \quad \IP{\gf,\del_\psi} =0 \mbox{ for all }\gf \in \wedge^\bullet T^*M.\]

So we propose that the relevant candidate higher Courant algebroid is the bundle:
\[C_\psi = TM \oplus (\wedge^{2n-2}T^*M \tensor \IP{\del_\psi}) \oplus \wedge^{2(n+k)-3}T^*M \oplus (\IP{\psi}\tensor \wedge^{2k-2}T^*M).\]
The logic behind the choices is that the first two summands (with vector field components) have degree $-1$ while the last two (which are just forms) have degree $2(n+k)-3$. A construction adding formal variables to $TM \oplus T^*M$ was studied before in the context of differential graded (dg) manifolds \cite{LUPERCIO201482}. Next we spell out the steps in non-dg language for higher Courant algebroids and without restricting our attention to specific numerical degrees of the forms involved.

Elements of $C_\psi$ act on $\wedge_\psi T^*M$ by interior product of (multi-)vectors and exterior product of forms:
\begin{equation}\label{eq:C_psi action}
\begin{aligned}
(X + \xi_{2n-2} \del_\psi + \psi \tensor \xi_{2k-2}+ \xi_{2(n+k)-3})& (\gf_0 + \psi \tensor \gf_1) : =\\ \iota_X \gf_0 + \xi_{2n-2}\wedge \gf_1 + \xi_{2(n+k)-3} \wedge \gf_0 &+\psi\tensor (-\iota_X \gf_1 -\xi_{2(n+k)-3} \wedge \gf_1+\xi_{2k-2}\wedge \gf_0).\\
\end{aligned}
\end{equation}
This acton gives rise to a natural pairing on $C_\psi$ using the graded commutator, $\{\cdot,\cdot\}$, which again corresponds to evaluation of (multi-)vectors on forms. To be precise, or $v\in C_\psi$ and $\gf \in \wedge^\bullet_\psi T^*M$,
\begin{equation}\label{eq:Clifford?}
\IP{v,v}\gf : = v \cdot(v \cdot\gf).
\end{equation}

Spelling it out: 
\begin{definition}\label{def:natural pairing}
The \emph{natural pairing} on $C_\psi$ is the fiberwise map
\[\IP{\cdot,\cdot} \colon C_\psi \times C_\psi \to \wedge^{2n-3} T^*M \tensor \IP{\del_\psi}\oplus \IP{\psi}\tensor\wedge^{2k-3} T^*M \oplus \wedge^{2(n+k)-4} T^*M,\]
\begin{equation}
\begin{aligned}2 \langle X + \xi_{2n-2} \del_\psi + \psi \tensor \xi_{2k-2}+& \xi_{2(n+k)-3}, Y + \eta_{2n-2} \del_\psi+ \psi \tensor \eta_{2k-2} + \eta_{2(n+k)-3}\rangle= \\
&=(\iota_X\eta_{2n-2} + \iota_Y \xi_{2n-2}) \del_\psi -\psi(\iota_X \eta_{2k-2} + \iota_Y \xi_{2k-2}) +\\
&+ \iota_X \eta_{2(n+k)-3} + \iota_Y \xi_{2(n+k)-3} + \xi_{2n-2}\eta_{2k-2} + \eta_{2n-2}\xi_{2k-2}.
\end{aligned}
\end{equation}
\end{definition}
Note that the mixed-term component in $\wedge^{2n-2}T^*M \tensor \IP{\del_\psi}$ pairs nontrivially with both vectors and the formal element $\psi$ that corresponds to the global angular form.

The action of $C_\psi$ on $\wedge^\bullet_\psi T^*M$ together with the twisted differential $d^H$ also gives rise to a Courant-like bracket:

\begin{theorem}[\cite{LUPERCIO201482}]
For a closed form $H\in \Omega^{2(n+k)-1}_\psi(M)$ and $v,w \in \Gamma(C_\psi)$, there is a unique element $\Cour{v,w}_H \in   \Gamma(C_\psi)$ such that following holds
\begin{equation}\label{eq:Cour1}
\Cour{v,w}_H\cdot \gf:= \{\{v, d^H\},w\}\gf\qquad \mbox{ for all } \gf\in\Omega_\psi(M),
\end{equation}
where $\{\cdot,\cdot\}$ denotes the graded commutator of operators.

Explicitly, for $H = H_{2(n+k)-1} + \psi\tensor H_{2k}$, $v = X + \xi_{2n-2} \del_\psi + \psi \tensor \xi_{2k-2}+ \xi_{2(n+k)-3}$, and $w = Y + \eta_{2n-2} \del_\psi + \psi \tensor \eta_{2k-2}+ \eta_{2(n+k)-3}$ we have
\begin{equation}\label{eq:Cour2}
\begin{aligned}
\Cour{v,w}_H &=    \iota_{[X,Y]}+(\L_{X}\eta_{2n-2} -\iota_{Y}(d\xi_{2n-2}+\iota_{X}\varepsilon)) \partial_\psi +\nonumber\\
&+\psi(\L_{X}\eta_{2k-2} -\iota_{Y}(d\xi_{2k-2}+\iota_{X}H_{2k})) +\nonumber \\
& +(\L_{X}\eta_{2(n+k)-3} -\iota_{Y}(d\xi_{2(n+k)-3}+\iota_{X}H_{2(n+k)-1}+\varepsilon \xi_{2k-2}+\xi_{2n-2}H_{2k})\nonumber\\
&+ (d\xi_{2n-2}+\iota_{X}\varepsilon)\eta_{2k-2} -(-1)^{2n-1}\eta_{2n-2} (d\xi_{2k-2}+\iota_{X}H_{2k}))\nonumber\\
\end{aligned}
\end{equation}
\end{theorem}


\begin{definition}\label{def:extended higher}
We call the triple $(C_\psi,\IP{\cdot,\cdot},\Cour{\cdot,\cdot}_{H})$ an \emph{extended Courant algebroid}.
\end{definition}
\begin{remark}
A more precise name would be \emph{extended higher Courant algebroid} which emphasises that the object being extended is a higher Courant algebroid. We will opt for the shorter name here.
\end{remark}

Much of the theory of Courant algebroids translates directly to extended Courant algebroids. For example, since $\Cour{\cdot,\cdot}_H$ is a derived bracket of odd degree, it satisfies a graded Jacobi identity, and hence sections of $C_\psi$ are derivations of this bracket:
\begin{equation}
\Cour{v_1,\Cour{v_2,v_3}_H}_H = \Cour{\Cour{v_1,v_2}_H,v_3}_H  + \Cour{v_2,\Cour{v_1,v_3}_H}_H.
\end{equation}
Similarly, it is a derivation of the natural pairing, in the sense that if we regard $\IP{v,w}$ as an operator that acts on $\wedge_\psi T^*M$ by \eqref{eq:Clifford?} and, for $v_1,v_2,v_3 \in \Gamma(C_\psi)$ define
\[\mc{L}^H_{v_1}\IP{v_2,v_3} := \frac{1}{2}\{\{v_1,d^H\},\{v_2,v_3\}\},\]
then it follows from the Jacobi identity for the graded commutator that
\begin{align*}
2\mc{L}^H_{v_1}\IP{v_2,v_3}  &:=\frac{1}{2} \{\{v_1,d^H\},\{v_2,v_3\}\}\\
&= \frac{1}{2}\{\{\{v_1,d^H\},v_2\},v_3\}+ \frac{1}{2}\{v_2,\{\{v_1,d^H\},v_3\}\}\\
&= \IP{\Cour{v_1,v_2}_H,v_3} +  \IP{v_2, \Cour{v_1,v_3}_H}.
\end{align*}
Which shows that sections of $C_\psi$ are also derivations of the natural pairing  and therefore are infinitesimal symmetries of the extended Courant algebroid.

The full study of the structure of the infinitesimal symmetries $C_\psi$ and related objects is beyond the scope of the present paper, except that $B$-field symmetries are relevant for T-duality.

\begin{definition}
An \emph{infinitesimal $B$-field transformation} of $C_\psi$ is the transformation of $C_\psi$ determined by a form $B\in \Omega_\psi^{2(n+k)-2}(M)$, given by $B(v) = \{B,v\}$.
\end{definition}
Explicitly, this translates to
\[B(X + \xi_{2n-2}\del_\psi + \psi \xi_{2k-2}+\xi_{2(n+k)-3}) = \psi \iota_XB_{2k-1} - \iota_XB_{2(n+k)-2} -\xi_{2n-2}B_{2k-1},\]
where $B = \psi B_{2k-1} + B_{2(n+k)-2}$.

Observe, however that the natural pairing is not invariant under the action of $B$-fields, but it is instead equivariant:
\[\IP{B(v),w)} - \IP{v,B(w)} = B(\IP{v,w}),\]
where for an element $\gf =  \gf_{2n-3} \del_\psi + \psi \gf_{2k-3} +  \gf_{2(n+k)-4}$ in the codomain of the natural pairing we have
\[B(\gf) = B_{2k-1}\gf_{2n-3}.\]

\subsection{Spherical T-duality as a map of extended Courant algebroids}

With the notion of extended Courant algebroid in place, we are ready to state and prove the theorem on geometric T-duality.

\begin{theorem}\label{theo:geometric}
Let $(E,H), (\hat{E},\hat{H}) \to M$ be a pair of spherical T-dual spaces with $dF =  p^*H -\hat{p}^*\hat{H}$ as in \eqref{eq:basic setup}. Let $\psi$ and $\hat\psi$ be global angular forms on $E$ and $\hat{E}$ and $C_\psi(E)$ and $C_{\hat\psi}(\hat{E})$ be the corresponding extended Courant algebroids.
Then 
\[\mc{T}_F\colon C_\psi(E) \to C_{\hat{\psi}}(\hat{E}), \qquad \mc{T}_F(v) = \hat{p}_*e^Fp^* v,\]
is a well defined map for which
\begin{equation}\label{eq:compatible}
\tau_F(v \cdot \gf) = \mc{T}_F(v) \cdot \tau_F(\gf), \qquad \mbox{ for all } v\in C_\psi, \gf \in \Omega_\psi.
\end{equation}
In particular
\begin{equation}\label{eq:courant isomorphim}
\mc{T}_F(\Cour{v,w}_H) = \Cour{\mc{T}_F(v),\mc{T}_F(w)}_{\hat{H}}.
\end{equation}
\end{theorem}
\begin{remark}
In the expression defining $\mc{T}_F$, the symbol $p^*v$ indicates all elements in $C_{\psi,\hat{\psi}}$ which are pull backs of $v$ via $p$, while $\hat{p}_*$ can only be applied to projectable elements of $C_{\psi,\hat{\psi}}$.
\end{remark}
\begin{proof}
By definition of the maps involved it follows that
\begin{itemize}
\item for $ w \in p^*v$, $p^*(v\cdot \gf) = w \cdot p^*(\gf)$,
\item $e^F(w \cdot \gf) =  e^Fw \cdot e^F(\gf)$
\item if $w\in C_{\psi,\hat{\psi}}$  is projectable under $\hat{p}_*$, then $\hat{p}_*(w\gf) = \hat{p}_*(w) \hat{p}_*(\gf)$.
\end{itemize}
Hence, if $\mc{T}_F$ is well defined, then \eqref{eq:compatible} follows.

To prove that at the level of Courant algebroids  the map $\hat{p}_* e^F p^*$ is well defined we will show that for each $v\in C_\psi$, there is precisely one element in $w \in p^*v$ for which $e^F w$ is projectable under $\hat{p}_*$ and therefore the assignment $v \mapsto \hat{p}_* w$ defines a bundle map.
 
Writing
\[F =  \hat{\psi}\wedge F_{2n-1} + \lambda \psi \wedge \hat{\psi} +\psi \wedge F_{2k-1} + F_{2(n+k)-2},\]
nondegeneracy of $F$ is equivalent to $\lambda \neq 0$. Given $v\in C_\psi$ if we fix an arbitrary element $w_0 \in p^*v$, any other element in $p^*v$ is of the form $w_0 + \eta_{2k-2}\del_{\hat{\psi}}$ and hence
\[e^F(w_0 + \eta_{2k-2}\del_{\hat{\psi}}) = e^F w_0 +\eta_{2k-2}\del_{\hat{\psi}} +\eta_{2k-2}(\lambda \psi - F_{2n-1}).\]
This is projectable,  if and only if it has no $\psi$ component, that is, if and only if taking interior product of the form part with $\del_\psi$ vanishes:
\[ \iota_{\del_\psi}(e^F w_0) + \lambda \eta_{2k-2} =0.\]
This holds if and only if we set $\eta_{2k-2}  = \frac{\iota_{\del_\psi}(e^F w_0)}{\lambda}$, showing that there is a unique $w \in p^*(v)$ with the desired property.

With equation \eqref{eq:compatible} at hand, \eqref{eq:courant isomorphim} is a formal algebraic consequence. Indeed, for all $\gf \in \Omega_\psi(E)$ we have
\begin{align*}
\tau_{F}(\Cour{v,w}_H \gf) &= \tau_F(\{\{v,d^H\},w\} \gf)\\
&= \{\{\mc{T}_Fv,d^{\hat{H}}\},\mc{T}_F{w}\}\tau_F(\gf)\\
& =\Cour{\mc{T}_Fv,\mc{T}_Fw}_{\hat{H}}\tau_F(\gf),
\end{align*}
where in the first equality we used the derived bracket description of $\Cour{\cdot,\cdot}_H$ and in the second equality we used Lemma \ref{lem:T-duality of subcomplexes} and equation \eqref{eq:compatible}.
Since $\tau_F$ is an isomorphism and the action of $C_\psi$ on $\Omega_\psi$ is faithful we conclude that \eqref{eq:courant isomorphim} holds.
\end{proof}

\begin{remark}
Since $B$-field transforms do not preserve the natural pairing (Definition \ref{def:natural pairing}), neither does T-duality. Instead, the pairing is equivariant due to the term $e^F$ in the definition of $\mc{T}_F$. 
\end{remark}

\begin{example} In this example we write down explicitly the map $\mc{T}_F$ for the T-dual pair constructed in Theorem \ref{theo:existence}. In that case we have, omitting pull backs, $F = \frac{1}{\lambda}\hat{\psi}\wedge \psi$. For this particular choice we have
\begin{equation}\label{eq:T-duality of Courant alg}
\mc{T}_F( X + \xi_{2n-2} \del_\psi + \psi \tensor \xi_{2k-2}+ \xi_{2(n+k)-3}) = X + \lambda \xi_{2k-2}\del_{\hat{\psi}} + \frac{1}{\lambda}\hat{\psi}\tensor \xi_{2n-2} + \xi_{2(n+k)-3}.
\end{equation}
Observe that even in this case, the natural pairing is not invariant under T-duality, unless we are in the unimodular case with $\lambda =1$.
\end{example}

\begin{example}[Exceptional Generalized Geometry] Here we recap an observation from \cite{LUPERCIO201482}, namely that some of the exceptional generalized geometries of \cite{Hull_2007,Pacheco_2008} reappear in the guise of extended Courant algebroids.

Consider $E = S^3 \times M \to M$ the trivial $S^3$-bundle with global angular form, $\psi$, the volume form on $S^3$ and twist the vanishing 7-form. Then the associated extended Courant algebroid over $M$ is
\[C_\psi =TM \oplus (\wedge^2 T^*M \tensor \IP{\del_\psi}) \oplus (\IP{\psi} \tensor \wedge^2 T^*M)  \oplus \wedge^5T^*M .\]
Since under the present hypothesis $E$ is self T-dual (following the construction of Theorem \ref{theo:existence}), and we can also consider the subspace of $C_\psi$ corresponding the self T-dual elements, that is,
\begin{align*}
C_\psi^{SD} &:= \{v \in C_\psi \colon \mc{T}_F v = v\}\\
&= \{X + \xi_2 \del_\psi + \psi \xi_2 + \xi_5\colon X \in TM, \xi_2 \in \wedge^2 T^*M, \xi_5 \in  \wedge^5T^*M\}
\end{align*}
Notice that since $\mc{T}_F$ preserves brackets, $C_\psi^{SD}$ is itself an involutive subbundle. Also  $\mathbb{T}M$ and $C_\psi^{SD}$ are isomorphic as bundles via the map
\[\mathbb{T}M \to C_{\psi}^{SD},\qquad  X + \xi_2 + \xi_5  \mapsto X + \sqrt{2}(\xi_2 \del_\psi + \psi \xi_2) + \xi_5.\]

This allows us to transport the bracket from $C_{\psi}$ to $\mathbb{T}M$ to obtain
\begin{align*}
\Cour{X + \xi_2 + \xi_5,Y + \eta_2+\eta_5} &=  [X,Y]+\left(\L_{X}\eta_{2} -\iota_{Y}(d\xi_{2})\right) +\left(\L_{X}\eta_{5} -\iota_{Y}d\xi_{5}+2(d\xi_{2})\eta_{2}\right),
\end{align*}
which agrees with the Dorfman bracket on $TM \oplus \wedge^2T^*M \oplus \wedge^5T^*M$ that appears in $E_6$ generalized geometry.
  
\end{example}

\section{Outlook}
There are a few points left untouched in the exposition above for which we do not have satisfactory answers.

\begin{question}\label{q:1} What is the mathematical relation between higher Courant algebroids and extended higher Courant algebroids?
\end{question}

From the mathematical point of view, given an oriented sphere bundle $E \to M$ it is still not clear what the relationship between the higher Courant algebroid $TE\oplus \wedge^{2(n+k)-3}T^*E$  defined on the total space of the sphere bundle and the extended higher Courant algebroid $C_\psi$ defined over $M$ is.

At first the method to obtain $C_\psi$ from $TE\oplus \wedge^{2(n+k)-3}T^*E$ is reminiscent of the method that produces heterotic Courant algebroids from  exact Courant algebroids in the context of principal bundles $P \to B$ \cite{MR3418511,MR3694654}. In the latter theory the relationship is made precise by stating that the heterotic Courant algebroid over $B$ is the result of a nonisotropic reduction of an exact Courant algebroid over $P$. For higher Courant algebroids it is not immediately clear that $C_\psi$ arises from $TE\oplus \wedge^{2(n+k)-3}T^*E$  by any reduction procedure.

\begin{question}
Is there a conceptual relation to exceptional Courant algebroids?
\end{question}

At the moment, the appearance of exceptional Courant algebroids as the self T-dual part of  extended higher Courant algebroid seems accidental. Notice that exceptional Courant algebroids come from representations of specific Lie groups and are restricted in which dimensions they may occur and the degrees of the forms one should consider, while extended higher Courant algebroids do not have those restrictions. Also the theory of exceptional Courant algebroids misses some of the symmetries and twists from our theory.

\begin{question}
Can string compactifications on spherical T-duals of different dimensions have related physical theories?
\end{question}

Most string theorists dislike\footnote{We heard from multiple sources that it makes no sense.} the idea of relating string theories on target spaces of different dimensions, yet, using extended Courant algebroids the relevant theory in both cases is encoded in bundles isomorphic to
\[TM \oplus \wedge^{2n-2}T^*M\oplus \wedge^{2k-2}T^*M \oplus \wedge^{2(n+k)-3}T^*M.\]
which contain few signs of the higher dimensional sphere bundle. From this point of view, the dimension of the sphere fiber might be ``hidden'' even for string theorists. The relation to $E_6$ generalized geometry suggests that this connection is not so farfetched.

\bibliographystyle{hyperamsplain-nodash}
\bibliography{references}

\providecommand{\bysame}{\leavevmode\hbox to3em{\hrulefill}\thinspace}
\providecommand{\MR}{\relax\ifhmode\unskip\space\fi MR }
\providecommand{\MRhref}[2]{%
  \href{http://www.ams.org/mathscinet-getitem?mr=#1}{#2}
}
\providecommand{\href}[2]{#2}
\begin{thebibliography}{10}

\bibitem{MR3418511}
D.~Baraglia and P.~Hekmati, \emph{Transitive {C}ourant algebroids, string
  structures and {$T$}-duality},
  \href{http://dx.doi.org/10.4310/ATMP.2015.v19.n3.a3}{Adv. Theor. Math. Phys.
  \textbf{19} (2015)}, no.~3, 613--672.

\bibitem{BV03}
I.~Belegradek and V.~Kapovitch, \emph{Obstructions to nonnegative curvature and
  rational homotopy theory},
  \href{http://dx.doi.org/10.1090/S0894-0347-02-00418-6}{J. Amer. Math. Soc.
  \textbf{16} (2003)}, no.~2, 259--284.

\bibitem{bott-tu}
R.~Bott and L.~W. Tu, \emph{Differential forms in algebraic topology}, Graduate
  Texts in Mathematics, vol.~82, Springer-Verlag, New York-Berlin, 1982.

\bibitem{BEM04}
P.~Bouwknegt, J.~Evslin, and V.~Mathai, \emph{{$T$}-duality: topology change
  from {$H$}-flux}, \href{http://dx.doi.org/10.1007/s00220-004-1115-6}{Comm.
  Math. Phys. \textbf{249} (2004)}, no.~2, 383--415.

\bibitem{BEM15}
P.~Bouwknegt, J.~Evslin, and V.~Mathai, \emph{Spherical {T}-duality},
  \href{http://dx.doi.org/10.1007/s00220-015-2354-4}{Comm. Math. Phys.
  \textbf{337} (2015)}, no.~2, 909--954.

\bibitem{BEM15b}
P.~Bouwknegt, J.~Evslin, and V.~Mathai, \emph{Spherical T-duality II: An
  infinity of spherical T-duals for non-principal SU(2)-bundles},
  \href{http://dx.doi.org/https://doi.org/10.1016/j.geomphys.2015.02.003}{Journal
  of Geometry and Physics \textbf{92} (2015)}, 46--54.

\bibitem{MR2246781}
U.~Bunke and T.~Schick, \emph{{$T$}-duality for non-free circle actions},
  Analysis, geometry and topology of elliptic operators, World Sci. Publ.,
  Hackensack, NJ, 2006, pp.~429--466.

\bibitem{Bus87}
T.~Buscher, \emph{A symmetry of the string background field equations},
  \href{http://dx.doi.org/https://doi.org/10.1016/0370-2693(87)90769-6}{Physics
  Letters B \textbf{194} (1987)}, no.~1, 59--62.

\bibitem{cavalcanti-2004}
G.~R. Cavalcanti, \emph{New aspects of the $dd^c$-lemma}, Ph.D. thesis, Oxford
  University, 2004.

\bibitem{CG11}
G.~R. Cavalcanti and M.~Gualtieri, \emph{Generalized complex geometry and
  T-duality}, A Celebration of the Mathematical Legacy of Raoul Bott (CRM
  Proceedings and Lecture Notes) (2010), 341--366.

\bibitem{CaWi22}
G.~R. Cavalcanti and A.~Witte, \emph{Non-principal T-duality, generalized
  complex geometry and blow-ups}, 2023.
  \href{http://arxiv.org/abs/2211.17173}{{\tt arXiv:2211.17173 [math.DG]}}.

\bibitem{Cortes:2021aab}
V.~Cort\'es and L.~David, \emph{{$T$-duality for transitive Courant
  algebroids}}, \href{http://dx.doi.org/10.4310/JSG.2023.v21.n4.a4}{J. Sympl.
  Geom. \textbf{21} (2023)}, no.~4, 775--856,
  \href{http://arxiv.org/abs/2101.07184}{{\tt arXiv:2101.07184 [math.DG]}}.

\bibitem{DELAOSSA1993377}
X.~C. {de la Ossa} and F.~Quevedo, \emph{Duality symmetries from non-abelian
  isometries in string theory},
  \href{http://dx.doi.org/https://doi.org/10.1016/0550-3213(93)90041-M}{Nuclear
  Physics B \textbf{403} (1993)}, no.~1, 377--394.

\bibitem{GARCIAFERNANDEZ20191059}
M.~Garcia-Fernandez, \emph{Ricci flow, Killing spinors, and T-duality in
  generalized geometry},
  \href{http://dx.doi.org/https://doi.org/10.1016/j.aim.2019.04.038}{Advances
  in Mathematics \textbf{350} (2019)}, 1059--1108.

\bibitem{MR3694654}
M.~Garcia-Fernandez, R.~Rubio, and C.~Tipler, \emph{Infinitesimal moduli for
  the {S}trominger system and {K}illing spinors in generalized geometry},
  \href{http://dx.doi.org/10.1007/s00208-016-1463-5}{Math. Ann. \textbf{369}
  (2017)}, no.~1-2, 539--595.

\bibitem{Hull_2007}
C.~M. Hull, \emph{Generalised geometry for M-theory},
  \href{http://dx.doi.org/10.1088/1126-6708/2007/07/079}{Journal of High Energy
  Physics \textbf{2007} (2007)}, no.~07, 079.

\bibitem{KlSe95}
C.~Klim\v{c}\'{\i}k and P.~\v{S}evera, \emph{Dual non-abelian duality and the
  {D}rinfel'd double},
  \href{http://dx.doi.org/10.1016/0370-2693(95)00451-P}{Phys. Lett. B
  \textbf{351} (1995)}, no.~4, 455--462.

\bibitem{LSW16}
J.~A. Lind, H.~Sati, and C.~Westerland, \emph{Twisted iterated algebraic
  {$K$}-theory and topological {T}-duality for sphere bundles},
  \href{http://dx.doi.org/10.2140/akt.2020.5.1}{Ann. K-Theory \textbf{5}
  (2020)}, no.~1, 1--42.

\bibitem{LUPERCIO201482}
E.~Lupercio, C.~Rengifo, and B.~Uribe, \emph{T-duality and exceptional
  generalized geometry through symmetries of dg-manifolds},
  \href{http://dx.doi.org/https://doi.org/10.1016/j.geomphys.2014.05.012}{Journal
  of Geometry and Physics \textbf{83} (2014)}, 82--98.

\bibitem{Pacheco_2008}
P.~P. Pacheco and D.~Waldram, \emph{M-theory, exceptional generalised geometry
  and superpotentials},
  \href{http://dx.doi.org/10.1088/1126-6708/2008/09/123}{Journal of High Energy
  Physics \textbf{2008} (2008)}, no.~09, 123.

\bibitem{RV92}
M.~Roček and E.~Verlinde, \emph{Duality, quotients, and currents},
  \href{http://dx.doi.org/https://doi.org/10.1016/0550-3213(92)90269-H}{Nuclear
  Physics B \textbf{373} (1992)}, no.~3, 630--646.

\end{thebibliography}

\end{document}